\documentclass[secnum,leqno,draft]{rims-bessatsu}
\usepackage{amssymb, amsmath}
\usepackage{graphics}
\usepackage[dvips]{graphicx}
\usepackage{eepic}
\usepackage{epic}
\newtheorem{theorem}{Theorem}[section]
\newtheorem{lemma}[theorem]{Lemma}
\newtheorem{proposition}[theorem]{Proposition}
\newtheorem{corollary}[theorem]{Corollary}

\newtheorem*{isom-prob}{Field isomorphism problem of a generic polynomial}
\newtheorem*{int-prob}{Field intersection problem of a generic polynomial}
\theoremstyle{definition}
\newtheorem*{definition}{Definition}
\newtheorem{example}[theorem]{Example}
\newtheorem{remark}[theorem]{Remark}
\newtheorem*{acknowledgments}{Acknowledgments}
\theoremstyle{remark}
\newcommand{\bs}{\mathbf{s}}\newcommand{\bt}{\mathbf{t}}
\newcommand{\ba}{\mathbf{a}}\newcommand{\bb}{\mathbf{b}}
\newcommand{\bx}{\mathbf{x}}\newcommand{\by}{\mathbf{y}}

\newcommand{\Gs}{G_\mathbf{s}}\newcommand{\Gt}{G_\mathbf{t}}
\newcommand{\Gst}{G_{\mathbf{s},\mathbf{t}}}
\newcommand{\Hst}{H_{\mathbf{s},\mathbf{t}}}

\newcommand{\opi}{\overline{\pi}}
\title{On the field intersection problem of generic polynomials: a survey}
\author{Akinari \textsc{Hoshi}\footnote{Department of Mathematics, Rikkyo University, 
3--34--1 Nishi-Ikebukuro Toshima-ku, Tokyo, 171-8501, Japan.}
~and Katsuya \textsc{Miyake}\footnote{Department of Mathematics, 
School of Fundamental Science and Engineering, Waseda University, 
3--4--1 Ohkubo Shinjuku-ku, Tokyo, 169-8555, Japan.}}
\AuthorHead{Akinari Hoshi and Katsuya Miyake}         
\support{This work was partially supported by Grant-in-Aid for Scientific 
Research (C) 19540057 of Japan Society for the Promotion of Science.}  
\VolumeNo{4x}           
\YearNo{200x}           
\PagesNo{000--000}      
\begin{document}
\maketitle
\begin{abstract}
Let $k$ be a field of characteristic $\neq 2$. 
We survey a general method of the field intersection problem of generic polynomials via 
formal Tschirnhausen transformation. 
We announce some of our recent results of cubic, quartic and quintic cases the details of 
which are to appear elsewhere. 
In this note, we give an explicit answer to the problem in the cases of cubic and dihedral 
quintic by using multi-resolvent polynomials. 
\end{abstract}
\section{Introduction}\label{seIntro}
Let $G$ be a finite group, $k$ a field of characteristic $\neq 2$, $M$ a field containing $k$ 
with $\# M=\infty$, and $k(\bt)$ the rational function field over $k$ with $n$ indeterminates 
$\bt=(t_1,\ldots,t_n)$. 
Our main interest in this note is a $k$-generic polynomial for $G$ 
(cf. \cite{DeM83}, \cite{Kem01}, \cite{JLY02}). 

\begin{definition}
A polynomial $f_\bt(X)\in k(\bt)[X]$ is called $k$-generic for $G$ if it has the following 
property: the Galois group of $f_\bt(X)$ over $k(\bt)$ is isomorphic to $G$ and every 
$G$-Galois extension $L/M$ over an arbitrary infinite field $M\supset k$ can be obtained as 
$L=\mathrm{Spl}_M f_\ba(X)$, the splitting field of $f_\ba(X)$ over $M$, for some 
$\ba=(a_1,\ldots,a_n)\in M^n$. 
\end{definition}

Let $f_\bt^G(X)\in k(\bt)[X]$ be a $k$-generic polynomial for $G$. 
Examples of $k$-generic polynomials for $G$ are known for various pairs of $(k,G)$ 
(for example, see \cite{Kem94}, \cite{KM00}, \cite{JLY02}, \cite{Rik04}). 
Since a $k$-generic polynomial $f_\bt^G(X)$ for $G$ covers all $G$-Galois extensions over 
$M\supset k$ by specializing parameters, it is natural to ask the following problem: 
%
\begin{isom-prob}
For a field $M\supset k$ and $\ba, \bb\in M^n$, 
determine whether $\mathrm{Spl}_M f_\ba^G(X)$ and $\mathrm{Spl}_M f_\bb^G(X)$ are 
isomorphic over $M$ or not. 
\end{isom-prob}
It would be desired to give an answer to the problem within the base field $M$ 
by using the data $\ba,\bb\in M^n$. 
Throughout this paper, we assume that $f_\ba^G(X)$ is separable for $\ba\in M^n$.

Let $S_n$ (resp. $D_n$, $C_n$) be the symmetric (resp. the dihedral, the cyclic) group 
of degree $n$. 
We take $k$-generic polynomials 
\begin{align*}
f_t^{C_3}(X)&:=X^3-tX^2-(t+3)X-1\,\in k(t)[X],\\
f_t^{S_3}(X)&:=X^3+tX+t\,\in k(t)[X],\\
f_{s,t}^{D_4}(X)&:=X^4+sX^2+t\,\in k(s,t)[X]
\end{align*}
for $G=C_3$, $S_3$, $D_4$, respectively. 
By using formal Tschirnhausen transformation, 
we showed the following theorem which is an analogue to 
the results of Morton \cite{Mor94} and Chapman \cite{Cha96}. 
\begin{theorem}[\cite{Mor94}, \cite{Cha96}, \cite{HM}]\label{thC3}
For $m,n\in M$, the splitting fields of $f_m^{C_3}(X)$ and 
of $f_n^{C_3}(X)$ over $M$ coincide if and only if there exists $z\in M$ such that either 
\[
n\,=\,\frac{m(z^3-3z-1)-9z(z+1)}{mz(z+1)+z^3+3z^2-1}\ \ \mathit{or}\ \ 
n\,=\,-\frac{m(z^3+3z^2-1)+3(z^3-3z-1)}{mz(z+1)+z^3+3z^2-1}.
\]
\end{theorem}
We also have analogues to the above theorem for two non-abelian groups $S_3$ and $D_4$ 
via formal Tschirnhausen transformation. 
\begin{theorem}[\cite{HM07}]\label{thS3isom}
Assume that char $k\neq 3$. 
For $a,b\in M$ with $a\neq b$, the splitting fields of $f_a^{S_3}(X)$ and of $f_b^{S_3}(X)$ 
over $M$ coincide if and only if there exists $u\in M$ such that 
\[
b=\frac{a(u^2+9u-3a)^3}{(u^3-2au^2-9au-2a^2-27a)^2}.
\]
\end{theorem}
\begin{theorem}[\cite{HM-2}]
For $a,b\in M$, we assume that $\mathrm{Gal}(f_{a,b}^{D_4}/M)=D_4$. 
Then for $a,b,a',b'\in M$, the splitting fields of $f_{a,b}^{D_4}(X)$ and of 
$f_{a',b'}^{D_4}(X)$ over $M$ coincide if and only if 
there exist $p,q\in M$ such that either 
\begin{align*}
{\rm (i)}\quad a'&=ap^2-4bpq+abq^2,\quad b'=b(p^2-apq+bq^2)^2\ \ \mathit{or}\\
{\rm (ii)}\quad a'&=2(ap^2-4bpq+abq^2),\quad b'=(a^2-4b)(p^2-bq^2)^2.
\end{align*}
\end{theorem}
By applying Hilbert's irreducibility theorem (cf. for example \cite[Chapter 3]{JLY02}) 
and Siegel's theorem for curves of genus $0$ (cf. for example \cite[Theorem 6.1]{Lan78}) 
to the theorems above respectively, we get the following corollaries: 
\begin{corollary}\label{cor0}
Let $f_\ba^G(X)=f_m^{C_3}(X)$ $($resp. $f_a^{S_3}(X)$, $f_{a,b}^{D_4}(X))$ be as above 
in Theorem {\rm 1.1} $($resp. Theorem {\rm 1.2}, Theorem {\rm 1.3}\,$)$ with given $\ba\in M^n$, 
and suppose that $M\supset k$ is Hilbertian $($e.g. a number field\,$)$. Then 
there exist infinitely many $\bb\in M^n$ such that 
$\mathrm{Spl}_{M} f_\ba^G(X)=\mathrm{Spl}_M f_\bb^G(X)$. 
\end{corollary}
\begin{corollary}\label{cor01}
Let $M$ be a number field and $\mathcal{O}_M$ the ring of integers in $M$. 
For $f_a^G(X)=f_m^{C_3}(X)$ $($resp. $f_a^{S_3}(X))$ as above in Theorem {\rm 1.1} 
$($resp. Theorem {\rm 1.2}\,$)$ with a given integer $a\in \mathcal{O}_M$, 
there exist only finitely many integers $b\in\mathcal{O}_M$ such that 
$\mathrm{Spl}_{M} f_a^G(X)=\mathrm{Spl}_M f_b^G(X)$. 
\end{corollary}
Indeed integers $b\in\mathcal{O}_M$ as in Corollary \ref{cor01} are 
derived from some integer solutions of (finitely many) cubic Thue equations 
$aXY(X+Y)+X^3+3XY^2-Y^3=\lambda$ (resp. $X^3-2aX^2Y-9aXY^2-2aY^3-27aY^3=\lambda$) 
which are curves of genus $1$ (see also the proof of \cite[Theorem 6.1]{Lan78}). 
\vspace*{2mm}

Kemper \cite{Kem01}, furthermore, showed that for a subgroup $H$ of $G$ every $H$-Galois 
extension over $M$ is also given by a specialization $f_\ba^G(X)$, $\ba \in M^n$, 
of a generic polynomial $f_\bt^G(X)$ for $G$. 
Hence a problem naturally arises. 
%
%
\begin{int-prob}
For a field $M\supset k$ and $\ba,\bb\in M^n$, determine 
the intersection of $\mathrm{Spl}_M f_\ba^G(X)$ and $\mathrm{Spl}_M f_\bb^G(X)$. 
\end{int-prob}
Clearly if we get an answer to the field intersection problem of a $k$-generic polynomial, 
we also obtain an answer to the corresponding field isomorphism problem. 

The aim of this note is to survey a method to give an answer to the isomorphism problem and 
the intersection problem of $k$-generic polynomials via formal Tschirnhausen transformation 
and multi-resolvent polynomials. 
In Section \ref{seResolv}, we review known results about resolvent polynomials. 
In Section \ref{seTschin}, we recall a formal Tschirnhausen transformation which is given 
in \cite{HM}. 
In Section \ref{seInt}, 
we give a general method to solve the intersection problem of $k$-generic polynomials. 
In Section \ref{seCubic}, we obtain an explicit answer to the problems in the cubic case. 
We give a proof of Theorem \ref{thS3isom} as a special case of the intersection problem 
of $f_t^{S_3}(X)=X^3+tX+t
$ via formal Tschirnhausen transformation. 
In Section \ref{seQuin}, we take the $k$-generic polynomial 
\begin{align*}
f_{s,t}^{D_5}(X)&=X^5+(t-3)X^4+(s-t+3)X^3+(t^2-t-2s-1)X^2+sX+t\in k(s,t)[X]
\end{align*}
for $D_5$ which is called Brumer's quintic (cf. \cite{JLY02}). 
Based on the general result, we illustrate how to give an answer to the problem 
for $f_{s,t}^{D_5}(X)$ by multi-resolvent polynomials. 
We also give some numerical examples. 

\section{Resolvent polynomials}\label{seResolv}

In this section we review some known results in the computational aspects of Galois theory 
(cf. the text books \cite{Coh93}, \cite{Ade01}). 
One of the fundamental tools in the determination of Galois group of a polynomial is resolvent 
polynomials; an absolute resolvent polynomial was introduced by Lagrange \cite{Lag1770} and 
a relative one by Stauduhar \cite{Sta73}. 
Several kinds of methods to compute resolvent polynomials have been developed by many 
mathematicians (see, for example, \cite{Sta73}, \cite{Gir83}, \cite{SM85}, \cite{Yok97}, 
\cite{MM97}, \cite{AV00}, \cite{GK00} and the references therein). 

Let $M(\supset k)$ be an infinite field and $\overline{M}$ a fixed algebraic closure of $M$.
Let $f(X):=\prod_{i=1}^m(X-\alpha_i) \in M[X]$ be a separable polynomial of degree $m$ with 
some fixed order of the roots $\alpha_1,\ldots,\alpha_m\in \overline{M}$. 
The Galois group of the splitting field $\mathrm{Spl}_M f(X)$ of $f(X)$ over $M$ 
may be obtained by using suitable resolvent polynomials. 

Let $k[\bx]:=k[x_1,\ldots,x_m]$ be the polynomial ring over $k$ with $m$ indeterminates 
$x_1,\ldots,x_m$. 
Put $R:=k[\bx, 1/\Delta_\bx]$, where $\Delta_\bx:=\prod_{1\leq i<j\leq m}(x_j-x_i)$. 
We take a surjective evaluation homomorphism $\omega_f : R \rightarrow 
k(\alpha_1,\ldots,\alpha_m),\, \Theta(x_1,\ldots,x_m)\mapsto 
\Theta(\alpha_1,\ldots,\alpha_m)$ for $\Theta \in R$.
We note that $\omega_f(\Delta_\bx)\neq 0$ from the assumption that $f(X)$ is separable. 
The kernel of the map $\omega_f$ is the ideal 
$I_f=\{\Theta(x_1,\ldots,x_m)\in R \mid \Theta(\alpha_1,\ldots,\alpha_m)=0\}$. 

Let $S_m$ be the symmetric group of degree $m$. 
We extend the action of $S_m$ on $m$ letters $\{1,\ldots,m\}$ 
to that on $R$ by $\pi(\Theta(x_1,\ldots,x_m)):=\Theta(x_{\pi(1)},\ldots,x_{\pi(m)})$. 
We define the Galois group of a polynomial $f(X)$ over $M$ by 
$\mathrm{Gal}(f/M):={\{\pi\in S_m \mid \pi(I_f)\subseteq I_f\}}$. 
We write $\mathrm{Gal}(f):=\mathrm{Gal}(f/M)$ for simplicity. 
The Galois group of the splitting field $\mathrm{Spl}_M f(X)$ of the polynomial $f(X)$ over $M$ 
is isomorphic to $\mathrm{Gal}(f)$. If we take another ordering of roots 
$\alpha_{\pi(1)},\ldots,\alpha_{\pi(m)}$ of $f(X)$ for some $\pi\in S_m$, the
corresponding realization of $\mathrm{Gal}(f)$ is conjugate in $S_m$. 
Hence, for arbitrary ordering of the roots of $f(X)$, $\mathrm{Gal}(f)$ 
is determined up to conjugacy in $S_m$. 

\begin{definition}
For $H\leq G\leq S_m$, an element $\Theta\in R$ is called a $G$-primitive $H$-invariant if 
$H=\mathrm{Stab}_G(\Theta)$ $:=$ $\{\pi\in G\ |\ \pi(\Theta)=\Theta\}$. 
For a $G$-primitive $H$-invariant $\Theta$, the polynomial 
\[
\mathcal{RP}_{\Theta,G}(X)=\prod_{\opi\in G/H}(X-\pi(\Theta))\in R^G[X], 
\]
where $\opi$ runs through the left cosets of $H$ in $G$, 
is called the {\it formal} $G$-relative $H$-invariant resolvent by $\Theta$, and the polynomial
\[
\mathcal{RP}_{\Theta,G,f}(X)=\prod_{\opi\in G/H}\bigl(X-\omega_f(\pi(\Theta))\bigr)\in M[X]
\]
is called the $G$-relative $H$-invariant resolvent of $f$ by $\Theta$. 
\end{definition}

The following theorem is fundamental in the theory of resolvent polynomials 
(cf. for example \cite[p.95]{Ade01}). 

\begin{theorem}\label{thfun}
Let $H\leq G\leq S_m$ be a tower of finite groups and $\Theta$ a $G$-primitive  $H$-invariant. 
Assume $\mathrm{Gal}(f)\leq G$. 
Suppose that $\mathcal{RP}_{\Theta,G,f}(X)=\prod_{i=1}^l h_i^{e_i}(X)$ gives the decomposition 
of $\mathcal{RP}_{\Theta,G,f}(X)$ into a product of powers of distinct irreducible 
polynomials $h_i(X), i=1,\ldots, l,$ in $M[X]$. Then we have a bijection 
\begin{align*}
\mathrm{Gal}(f)\backslash G/H\quad &\longrightarrow \quad \{h_1^{e_1}(X),\ldots,h_l^{e_l}(X)\}\\
\mathrm{Gal}(f)\, \pi\, H\quad &\longmapsto\quad h_\pi(X)
=\prod_{\tau H\subseteq \mathrm{Gal}(f)\,\pi\,H}\bigl(X-\omega_{f}(\tau(\Theta))\bigr)
\end{align*}
where the product runs through the left cosets $\tau H$ of $H$ in $G$ contained in 
$\mathrm{Gal}(f)\, \pi\, H$, that is, through $\tau=\pi_\sigma \pi$ where $\pi_\sigma$ runs 
a system of representative of the left cosets of $\mathrm{Gal}(f) \cap \pi H\pi^{-1};$ each 
$h_\pi(X)$ is irreducible or a power of an irreducible polynomial with $\mathrm{deg}(h_\pi(X))$ 
$=$ $|\mathrm{Gal}(f)\, \pi\, H|/|H|$ $=$ $|\mathrm{Gal}(f)|/|\mathrm{Gal}(f)\cap \pi H\pi^{-1}|$. 
\end{theorem}

\begin{corollary} 
If $\mathrm{Gal}(f)\leq \pi H\pi^{-1}$ for some $\pi\in G$ then 
$\mathcal{RP}_{\Theta,G,f}(X)$ has a linear factor over $M$. 
Conversely, if $\mathcal{RP}_{\Theta,G,f}(X)$ has a non-repeated linear factor over $M$ 
then there exists $\pi\in G$ such that $\mathrm{Gal}(f)\leq \pi H\pi^{-1}$. 
\end{corollary}

\begin{remark}\label{remGir}
When the resolvent polynomial $\mathcal{RP}_{\Theta,G,f}(X)$ has a repeated factor, there 
always exists a suitable Tschirnhausen transformation $\hat{f}$ of $f$ (cf. \S 3) over $M$ 
(resp. $X-\hat{\Theta}$ of $X-\Theta$ over $k$) such that $\mathcal{RP}_{\Theta,G,\hat{f}}(X)$ 
(resp. $\mathcal{RP}_{\hat{\Theta},G,f}(X)$) has no repeated factors (cf. \cite{Gir83}, 
\cite[Alg. 6.3.4]{Coh93}, \cite{Col95}). 
\end{remark}

In the case where $\mathcal{RP}_{\Theta,G,f}(X)$ has no repeated factors, 
we have the followings: \\
{\rm (i)} For $\pi\in G$, the fixed group of the field $M\bigl(\omega_{f}(\pi(\Theta))\bigr)$ 
corresponds to $\mathrm{Gal}(f)\cap \pi H\pi^{-1}$. 
In particular, we have $\mathrm{Gal}(\mathcal{RP}_{\Theta,G,f})\cong \mathrm{Gal}(f)/N$ with 
$N=\mathrm{Gal}(f)\cap \bigcap_{\pi\in G}\pi H\pi^{-1}$; \\
{\rm (ii)} let $\varphi : G\rightarrow S_{[G:H]}$ denote the permutation representation of 
$G$ on the set of left cosets $G/H$ given by the left multiplication. Then we 
have a realization of the Galois group of $\mathrm{Spl}_M \mathcal{RP}_{\Theta,G,f}(X)$ 
as a subgroup of the symmetric group $S_{[G:H]}$ by $\varphi(\mathrm{Gal}(f))$. 

\section{Formal Tschirnhausen transformation}\label{seTschin}

We recall the geometric interpretation of a Tschirnhausen transformation 
which is given in \cite{HM} (see also \cite{HM-3}). 
Let $f(X)$ be a monic separable polynomial of degree $n$ in $M[X]$ 
with a fixed order of the roots $\alpha_1,\ldots,\alpha_n$ of $f(X)$ in $\overline{M}$. 
A Tschirnhausen transformation of $f(X)$ over $M$ is a polynomial of the form 
\[
g(X)=\prod_{i=1}^n 
\bigl(X-(c_0+c_1\alpha_i+\cdots+c_{n-1}\alpha_i^{n-1})\bigr),\ c_i \in M.
\]
Two polynomials $f(X)$ and $g(X)$ in $M[X]$ are Tschirnhausen equivalent over $M$ if they are 
Tschirnhausen transformations over $M$ of each other. 
For two irreducible separable polynomials $f(X)$ and $g(X)$ in $M[X]$, 
$f(X)$ and $g(X)$ are Tschirnhausen equivalent over $M$ if and only if 
the quotient fields $M[X]/(f(X))$ and $M[X]/(g(X))$ are $M$-isomorphic. 

In order to obtain an answer to the field intersection problem of $k$-generic polynomials 
via multi-resolvent polynomials, we first treat a general polynomial 
whose roots are $n$ indeterminates $x_1,\ldots,x_n$: 
\[
f_\bs(X)\, =\, \prod_{i=1}^n(X-x_i)\, 
=\, X^n-s_1X^{n-1}+s_2X^{n-2}+\cdots+(-1)^n s_n\ \in k[\bs][X]
\]
where $k[x_1,\ldots,x_n]^{S_n}=k[\bs]:=k[s_1,\ldots,s_n], \bs=(s_1,\ldots,s_n),$ and $s_i$ is the 
$i$-th elementary symmetric function in $n$ variables $\bx=(x_1,\ldots,x_n)$. 

Put $R_{\bx,\by}:=k[\bx,\by,1/\Delta_\bx,1/\Delta_\by]$, where $\by=(y_1,\ldots,y_n)$ consists of 
$n$ indeterminates, 
$\Delta_\bx:=\prod_{1\leq i<j\leq n}(x_j-x_i)$ and $\Delta_\by:=\prod_{1\leq i<j\leq n}(y_j-y_i)$. 
We define the interchanging involution $\iota_{\bx,\by}$ which exchanges the 
$x_i$'s and the $y_i$'s: 
\begin{align}
\iota_{\bx,\by}\ :\ R_{\bx,\by}\longrightarrow R_{\bx,\by},\ 
x_i\longmapsto y_i,\ y_i\longmapsto x_i,\quad (i=1,\ldots,n).\label{defiota}
\end{align}
We take another general polynomial $f_\bt(X):=\iota_{\bx,\by}(f_\bs(X))\in k[\bt][X], 
\bt=(t_1,\ldots, t_n)$, whose roots are $n$ indeterminates $y_1,\ldots,y_n$ where 
$t_i=\iota_{\bx,\by}(s_i)$ is the $i$-th elementary symmetric function in $\by=(y_1,\ldots,y_n)$. 
We put $K:=k(\bs,\bt)$ and $f_{\bs,\bt}(X):=f_\bs(X)f_\bt(X)$. 
The polynomial $f_{\bs,\bt}(X)$ of degree $2n$ is defined over $K$. We denote 
\[
\Gs\, :=\, \mathrm{Gal}(f_\bs/K),\quad \Gt\, :=\, \mathrm{Gal}(f_\bt/K),\quad 
\Gst\, :=\, \mathrm{Gal}(f_{\bs,\bt}/K). 
\]
Then we have $\Gst=\Gs\times\Gt, \Gs\cong \Gt\cong S_n$ and $k(\bx,\by)^{\Gst}=K$. 

We intend to apply the results of the previous section for 
$m=2n$, $G=\Gst\leq S_{2n}$ and $f=f_{\bs,\bt}$. 

Note that in the splitting field $\mathrm{Spl}_K f_{\bs,\bt}(X)=k(\bx,\by)$, there exist $n!$ 
Tschirnhausen transformations from $f_\bs(X)$ to $f_\bt(X)$ with respect to 
$y_{\pi(1)},\ldots,y_{\pi(n)}$ for $\pi\in S_n$. 
We shall study the field of definition of each Tschirnhausen transformation from 
$f_\bs(X)$ to $f_\bt(X)$. 
Let $D:=[x_i^{j-1}]_{1\leq i,j\leq n}$ be the Vandermonde matrix of size $n$. 
The matrix $D\in M_n(k(\bx))$ is invertible because ${\rm det}\,D=\Delta_\bx$. 
The field $k(\bs)(\Delta_\bx)$ is a quadratic extension of $k(\bs)$ 
which corresponds to the fixed field of the alternating group of degree $n$. 
We define the $n$-tuple $(u_0(\mathbf{x},\mathbf{y}),\ldots,
u_{n-1}(\mathbf{x},\mathbf{y}))\in (R_{\bx,\by})^n$ by 
\[
\left[\begin{array}{c}u_0(\mathbf{x},\mathbf{y})\\ u_1(\mathbf{x},
\mathbf{y})\\ \vdots \\ u_{n-1}(\mathbf{x},\mathbf{y})\end{array}\right]
:=D^{-1}\left[\begin{array}{c}y_1\\ y_2\\ \vdots \\ 
y_n\end{array}\right]. 
\]
Cramer's rule shows us 
\[
u_i(\mathbf{x},\mathbf{y})=\Delta_\bx^{-1}\cdot\mathrm{det}
\left[\begin{array}{cccccccc}
1 & x_1 & \cdots & x_1^{i-1} & y_1 & x_1^{i+1} & \cdots & x_1^{n-1}\\ 
1 & x_2 & \cdots & x_2^{i-1} & y_2 & x_2^{i+1} & \cdots & x_2^{n-1}\\
\vdots & \vdots & & \vdots & \vdots & \vdots & & \vdots\\
1 & x_n & \cdots & x_n^{i-1} & y_n & x_n^{i+1} & \cdots & x_n^{n-1}
\end{array}\right].
\]
We write $u_i:=u_i(\mathbf{x},\mathbf{y}), (i=0,\ldots,n-1)$. 
The Galois group $\Gst$ acts on the orbit $\{\pi(u_i)\ |\ \pi\in \Gst \}$ 
via regular representation from the left. 
However this action is not faithful. 
We put 
\[
\Hst:=\{(\pi_\bx, \pi_\by)\in \Gst\ |\ \pi_\bx(i)=\pi_\by(i)\ \mathrm{for}\ 
i=1,\ldots,n \}\cong S_n. 
\]
If $\pi \in \Hst$ then we have $\pi(u_i)=u_i$ for $i=0,\ldots,n-1$. 
Indeed we see the following: 
\begin{lemma}\label{stabil}
For $i$, $0\leq i\leq n-1$, $u_i$ is a $\Gst$-primitive $\Hst$-invariant. 
\end{lemma}

Let $\Theta:=\Theta(\bx,\by)$ be a $\Gst$-primitive $\Hst$-invariant. 
Let $\opi=\pi\Hst$ be a left coset of $\Hst$ in $\Gst$. 
The group $\Gst$ acts on the set $\{ \pi(\Theta)\ |\ \opi\in \Gst/\Hst\}$ 
transitively from the left through the action on the set $\Gst/\Hst$ of left cosets. 
Each of the sets $\{ \overline{(1,\pi_\by)}\ |\ (1,\pi_\by)\in \Gst\}$ 
and $\{ \overline{(\pi_\bx,1)}\ |\ (\pi_\bx,1)\in \Gst\}$ forms a complete residue 
system of $\Gst/\Hst$, and hence the subgroups $\Gs$ and $\Gt$ of $\Gst$ act on the set 
$\{ \pi(\Theta)\ |\ \opi\in \Gst/\Hst\}$ transitively. 
For $\opi=\overline{(1,\pi_\by)}\in \Gst/\Hst$, we obtain the following equality: 
\[
y_{\pi_\by(i)} = \pi_\by(u_0)+\pi_\by(u_1) x_i+\cdots+\pi_\by(u_{n-1})x_i^{n-1}\ 
\mathrm{for}\ i=1,\ldots,n. 
\]
Hence the set $\{(\pi(u_0),\ldots,\pi(u_{n-1}))\ |\ \opi\in \Gst/\Hst\}$ 
gives coefficients of $n!$ different Tschirnhausen transformations from $f_\bs(X)$ 
to $f_\bt(X)$ each of which is defined over $K(\pi(u_0),\ldots,\pi(u_{n-1}))$ respectively. 
\begin{definition}
We call $K(\pi(u_0),\ldots,\pi(u_{n-1})), (\opi\in \Gst/\Hst)$, a field of formal 
Tschirnhausen coefficients from $f_\bs(X)$ to $f_\bt(X)$. 
\end{definition}
Put $v_i:=\iota_{\bx,\by}(u_i)$, for $i=0,\ldots,n-1$. 
Then $v_i$ is also a $\Gst$-primitive $\Hst$-invariant, and $K(\pi(v_0),\ldots,\pi(v_{n-1}))$ 
gives a field of formal Tschirnhausen coefficients from $f_\bt(X)$ to $f_\bs(X)$. 
\begin{proposition}\label{prop1}
For every  $\Gst$-primitive $\Hst$-invariant $\Theta$, we have $K(\pi(\Theta))$ $=$ 
$k(\bx,\by)^{\pi\Hst \pi^{-1}}\!\!$ $=$ $K(\pi(u_0),\ldots,\pi(u_{n-1}))$, 
and $[K(\pi(\Theta)) : K]=n!$ for each $\opi\in \Gst/\Hst$. 
\end{proposition}
Hence, for each of $n!$ fields $K(\pi(\Theta))$, 
we have $\mathrm{Spl}_{K(\pi(\Theta))} f_\bs(X)=\mathrm{Spl}_{K(\pi(\Theta))} f_\bt(X)$. 
\begin{proposition}\label{propLL}
Let $\Theta$ be a $\Gst$-primitive $\Hst$-invariant. 
Then we have 
\begin{align*}
&{\rm (i)}\ K(\bx)\cap K(\pi(\Theta))=K(\by)\cap K(\pi(\Theta))=K\quad \textrm{for}\quad 
\opi\in \Gst/\Hst\,{\rm ;}\\
&{\rm (ii)}\ K(\bx,\by)=K(\bx,\pi(\Theta))=K(\by,\pi(\Theta))\quad \textrm{for}\quad 
\opi\in \Gst/\Hst\,{\rm ;}\\
&{\rm (iii)}\ K(\bx,\by)=K(\pi(\Theta)\ |\ \opi\in \Gst/\Hst). 
\end{align*}
\end{proposition}

\section{Field intersection problem}\label{seInt}

In this section, we explain how to get an answer to the field intersection problem of 
generic polynomials via multi-resolvent polynomial (cf. \cite{HM}, \cite{HM-3}). 
For $\ba=(a_1,\ldots,a_n)$, $\bb=(b_1,\ldots,b_n)\in M^n$, 
we take some fixed order of the roots $\alpha_1,\ldots,\alpha_n$ of $f_\ba(X)$ 
and $\beta_1,\ldots,\beta_n$ of $f_\bb(X)$ in $\overline{M}$, respectively, and denote 
$L_{\ba} :=M(\alpha_1,\ldots,\alpha_n)$ and $L_{\bb}:=M(\beta_1,\ldots,\beta_n)$. 
We put $f_{\ba,\bb}(X):=f_\ba(X)f_\bb(X)\in M[X]$ and define a specialization homomorphism 
$\omega_{f_{\ba,\bb}}$ by 
\[
\omega_{f_{\ba,\bb}}\, :\ R_{\bx,\by}\,\longrightarrow\,
k(\alpha_1,\ldots,\alpha_n,\beta_1,\ldots,\beta_n),\quad 
\Theta(\bx,\by)\,\longmapsto\,\Theta(\alpha_1,\ldots,\alpha_n,\beta_1,\ldots,\beta_n). 
\]
Put $\Delta_\ba:=\omega_{f_{\ba,\bb}}(\Delta_\bx)$ and 
$\Delta_\bb:=\omega_{f_{\ba,\bb}}(\Delta_\by)$. 
We assume that both of the polynomials $f_\ba(X)$ and $f_\bb(X)$ are separable over $M$, 
i.e. $\Delta_\ba\cdot\Delta_\bb\neq 0$. 
Put $G_\ba:=\mathrm{Gal}(f_{\ba}/M)$, $G_\bb:=\mathrm{Gal}(f_{\bb}/M)$ and 
$G_{\ba,\bb}:=\mathrm{Gal}(f_{\ba,\bb}/M)$. 
Then we may naturally regard $G_{\ba,\bb}$ as a subgroup of $\Gst$. 
For $\opi \in \Gst/\Hst$, we put $c_{j,\pi}:=\omega_{f_{\ba,\bb}}(\pi(u_j)), 
d_{j,\pi}:=\omega_{f_{\ba,\bb}}\bigl(\pi(\iota_{\bx,\by}(u_j))\bigr)$, ($j=0,\ldots,n-1$). 
Then for each $i = 1, \ldots, n$, we have 
\begin{align*}
\beta_{\pi_\by(i)}\,&=\, c_{0,\pi} + c_{1,\pi}\,\alpha_{\pi_\bx(i)}+ \cdots 
+ c_{n-1,\pi}\,\alpha_{\pi_\bx(i)}^{n-1}, \\
\alpha_{\pi_\bx(i)}\,&=\, d_{0,\pi} + d_{1,\pi}\,\beta_{\pi_\by(i)}+ \cdots 
+ d_{n-1,\pi}\,\beta_{\pi_\by(i)}^{n-1}.
\end{align*}
For each $\opi\in\Gst/\Hst$, there exists a Tschirnhausen transformation from $f_\ba(X)$ 
to $f_\bb(X)$ over the field of Tschirnhausen coefficients $M(c_{0,\pi},\ldots,c_{n-1,\pi})$. 
From the assumption $\Delta_\ba\cdot \Delta_\bb\neq 0$, we first see the following lemma 
(cf. \cite[p. 141]{JLY02}, \cite{HM}). 
\begin{lemma}\label{lemM}
Let $M'/M$ be a field extension. 
If $f_\bb(X)$ is a Tschirnhausen transformation of $f_\ba(X)$ over $M'$, then $f_\ba(X)$ 
is a Tschirnhausen transformation of $f_\bb(X)$ over $M'$. 
In particular, we have $M(c_{0,\pi},\ldots,c_{n-1,\pi})=M(d_{0,\pi},\ldots,d_{n-1,\pi})$ 
for every $\opi\in\Gst/\Hst$. 
\end{lemma}

To obtain an answer to the field intersection problem of $f_\bs(X)$ we study the 
$n!$ fields 
$M(c_{0,\pi},\ldots,c_{n-1,\pi})$ of Tschirnhausen coefficients from $f_\ba(X)$ to 
$f_\bb(X)$ over $M$. 
\begin{proposition}\label{propc}
Under the assumption $\Delta_\ba\cdot\Delta_\bb\neq 0$, we have\\
{\rm (i)}\ \ $\mathrm{Spl}_{M(c_{0,\pi},\ldots,c_{n-1,\pi})} f_\ba(X)
=\mathrm{Spl}_{M(c_{0,\pi},\ldots,c_{n-1,\pi})} f_\bb(X)\ \,\textit{for each}\ \,
\opi\in\Gst/\Hst$\,{\rm ;}\\
{\rm (ii)}\ $L_\ba L_\bb=L_\ba\, M(c_{0,\pi},\ldots,c_{n-1,\pi})=L_\bb\, M(c_{0,\pi},\ldots,c_{n-1,\pi})$\ \,for each\ \,$\opi\in\Gst/\Hst$.
\end{proposition}
Applying the specialization $\omega_{f_{\ba,\bb}}$, 
we take the $\Gst$-relative $\Hst$-invariant resolvent polynomial of $f_{\ba,\bb}$ 
by a $\Gst$-primitive $\Hst$-invariant $\Theta$: 
\[
\mathcal{RP}_{\Theta,\Gst,f_{\ba,\bb}}(X)\, =\, 
\prod_{\opi\in \Gst/\Hst} \bigl(X-\omega_{f_{\ba,\bb}}(\pi(\Theta))\bigr)\in M[X]. 
\]
The resolvent polynomial $\mathcal{RP}_{\Theta,\Gst,f_{\ba,\bb}}(X)$ is also called 
an (absolute) multi-resolvent (cf. \cite{GLV88}, \cite{RV99}, \cite{Val}). 
\begin{proposition}\label{prop12}
For $\ba,\bb \in M^n$ with $\Delta_\ba\cdot\Delta_\bb\neq 0$, suppose that the polynomial 
$\mathcal{RP}_{\Theta,\Gst,f_{\ba,\bb}}(X)$ has no repeated factors. 
Then the following two assertions hold\,{\rm :}\\
$(\mathrm{i})$\ $M(c_{0,\pi},\ldots,c_{n-1,\pi})=M\bigl(\omega_{f_{\ba,\bb}}(\pi(\Theta))\bigr)$\, 
for each $\,\opi\in\Gst/\Hst$\,{\rm ;}\\
$(\mathrm{ii})$\ ${\rm Spl}_M f_{\ba,\bb}(X)
=M(\omega_{f_{\ba,\bb}}(\pi(\Theta))\ |\ \opi\in\Gst/\Hst)$. 
\end{proposition}

\begin{definition}
For a separable polynomial $f(X)\in M[X]$ of degree $d$, the decomposition type of $f(X)$ 
over $M$, denoted by {\rm DT}$(f/M)$, is defined as the partition of $d$ induced by the 
degrees of the irreducible factors of $f(X)$ over $M$. 
We define the decomposition type {\rm DT}$(\mathcal{RP}_{\Theta,G,f}/M)$ of 
$\mathcal{RP}_{\Theta,G,f}(X)$ over $M$ by {\rm DT}$(\mathcal{RP}_{\Theta,G,\hat{f}}/M)$ where 
$\hat{f}(X)$ is a Tschirnhausen transformation of $f(X)$ over $M$ which satisfies that 
$\mathcal{RP}_{\Theta,G,\hat{f}}(X)$ has no repeated factors (cf. Remark \ref{remGir}). 
\end{definition}

We write $\mathrm{DT}(f):=\mathrm{DT}(f/M)$ for simplicity. 
From Theorem \ref{thfun}, the decomposition type 
$\mathrm{DT}(\mathcal{RP}_{\Theta,\Gst,f_{\ba,\bb}})$ coincides with the partition of $n!$ 
induced by the lengths of the orbits of $\Gst/\Hst$ under the action of 
$\mathrm{Gal}(f_{\ba,\bb})$. 
Hence, by Proposition \ref{prop12}, $\mathrm{DT}(\mathcal{RP}_{\Theta,\Gst,f_{\ba,\bb}})$ 
gives the degrees of $n!$ fields of Tschirnhausen coefficients 
from $f_\ba(X)$ to $f_\bb(X)$ over $M$.  

We conclude that $\mathrm{DT}(\mathcal{RP}_{\Theta,\Gst,f_{\ba,\bb}})$ gives us information 
about the field intersection problem for $f_\bs(X)$ through the degrees of the fields 
of Tschirnhausen coefficients $M(c_{0,\pi},\ldots,c_{n-1,\pi})$ over $M$ 
and is determined by the degeneration of the Galois group $\mathrm{Gal}(f_{\ba,\bb})$ 
under the specialization $(\bs, \bt) \mapsto (\ba, \bb)$. 
As a special case of the field intersection problem, we get the followings: 
\begin{theorem}[\cite{HM-3}]\label{throotf}
For $\ba,\bb \in M^n$ with $\Delta_\ba\cdot\Delta_\bb\neq 0$, 
the quotient fields $M[X]/(f_\ba(X))$ and $M[X]/(f_\bb(X))$ are $M$-isomorphic if and only if 
the decomposition type ${\rm  DT}(\mathcal{RP}_{\Theta,\Gst,f_{\ba,\bb}})$ over $M$ includes $1$.
\end{theorem}
\begin{corollary}[The field isomorphism problem]\label{cor1}
For $\ba,\bb \in M^n$ with $\Delta_\ba\cdot\Delta_\bb\neq 0$, we assume that 
both of $G_\ba$ and $G_\bb$ are isomorphic to a transitive subgroup $G\leq S_n$ 
and that all subgroups of $G$ with index $n$ are conjugate in $G$. 
Then $\mathrm{DT}(\mathcal{RP}_{\Theta,\Gst,f_{\ba,\bb}})$ over $M$ includes $1$ if and only if 
$\mathrm{Spl}_M f_\ba(X)$ and $\mathrm{Spl}_M f_\bb(X)$ coincide. 
\end{corollary}

\section{Field intersection problem: the case of $S_3$}\label{seCubic}

We take $f_s^{S_3}(X)=X^3+sX+s\in k(s)[X]$. 
For $a,b\in M$, we put $L_a:=\mathrm{Spl}_M f_a^{S_3}(X)$ and 
$G_a:=\mathrm{Gal}(f_a^{S_3}/M), G_{a,b}:=\mathrm{Gal}(f_a^{S_3}f_b^{S_3}/M)$. 
For $a,b\in M$, we assume that $ab(4a+27)(4b+27)\neq 0$ because the discriminant of 
$f_s^{S_3}(X)$ equals $-s^2(4s+27)$. 

In the case of char $k\neq 3$, we take a $\Gst$-primitive $\Hst$-invariant 
$\Theta:=3u_1/u_2$ where 
$u_1$ and $u_2$ are formal Tschirnhausen coefficients which are given in Section \ref{seTschin}. 
Then we may evaluate 
\begin{align*}
F_{s,t}(X)&:=(s-t)\cdot\mathcal{RP}_{\Theta,\Gst,f_s^{S_3}f_t^{S_3}}(X)
=(s-t)\cdot\!\!\prod_{\opi\in\Gst/\Hst}\!\!(X-\pi(\Theta))\, \in k(s,t)[X]\\
&\,=\,s(X^2+9X-3s)^3-t(X^3-2sX^2-9sX-2s^2-27s)^2. 
\end{align*}
The discriminant of $F_{s,t}(X)$ is $s^{10}t^4(4s+27)^{15}(4t+27)^3$. 
Note that, for $a,b\in M$ with $ab(4a+27)(4b+27)\neq 0$, 
$F_{a,b}(X)$ has no repeated factors. 
In the case of char $k=3$, we may take $\Theta=u_0$, 
the constant term of a formal Tschirnhausen transformation, as a suitable $\Gst$-primitive 
$\Hst$-invariant. 
Then we get squarefree $\mathcal{RP}_{u_0,\Gst,f_s^{S_3}f_t^{S_3}}(X)$ 
as we mentioned in Remark \ref{remGir} (see \cite{HM}). 
\begin{theorem}[\cite{HM}]\label{thS3}
Assume that char $k\neq 3$. 
For $a,b\in M$, we also assume that $a\neq b$ and $\#G_a\geq \#G_b>1$. 
Then an answer to the field intersection problem for $f_s^{S_3}(X)$ 
is given by {\rm DT}$(F_{a,b})$ as Table $1$ shows. 
\end{theorem}
\begin{center}
{\small 
{\rm Table} $1$\vspace*{2mm}\\
}
{\small 
\begin{tabular}{|c|c|c|l|l|}\hline
$G_a$& $G_b$ & $G_{a,b}$ & & ${\rm DT}(F_{a,b})$\\ \hline 
& & $S_3\times S_3$ & $L_a\cap L_b=M$ & $6$ \\ \cline{3-5} 
& $S_3$ & $(C_3\times C_3)\rtimes C_2$ & $[L_a\cap L_b : M]=2$ & $3,3$\\ \cline{3-5}
\raisebox{-1.6ex}[0cm][0cm]{$S_3$} & & $S_3$ & $L_a=L_b$ & $3,2,1$ \\ \cline{2-5}
& $C_3$ & $S_3\times C_3$ & $L_a\cap L_b=M$ & $6$ \\ \cline{2-5}
& \raisebox{-1.6ex}[0cm][0cm]{$C_2$} & $D_6$ & $L_a\not\supset L_b$ & $6$\\ \cline{3-5}
& & $S_3$ & $L_a\supset L_b$ & $3,3$ \\ \hline
& \raisebox{-1.6ex}[0cm][0cm]{$C_3$} & $C_3\times C_3$ & $L_a\neq L_b$ & $3,3$ \\ \cline{3-5}
$C_3$ & & $C_3$ & $L_a=L_b$ & $3,1,1,1$\\ \cline{2-5}
& $C_2$ & $C_6$ & $L_a\cap L_b=M$ & $6$ \\ \hline
\raisebox{-1.6ex}[0cm][0cm]{$C_2$} & \raisebox{-1.6ex}[0cm][0cm]{$C_2$} & $C_2\times C_2$ 
& $L_a\neq L_b$ & $4,2$ \\ \cline{3-5}
& & $C_2$ & $L_a=L_b$ & $2,2,1,1$ \\ \hline
\end{tabular}
}\vspace*{1mm}
\end{center}
As a special case of Theorem \ref{thS3} (cf. Theorem \ref{throotf} and Corollary 
\ref{cor1}), we get Theorem \ref{thS3isom} which we introduced in Section \ref{seIntro}. 
\vspace*{2mm}

{\bf Proof of Theorem \ref{thS3isom}.}
Because the polynomial 
$
F_{s,t}(X)=s(X^2+9X-3s)^3-t(X^3-2sX^2-9sX-2s^2-27s)^2
$
is linear in $t$, we see that $F_{a,b}(X)$ has a root in $M$ if and only if there exists 
$u\in M$ such that
\[
b=\frac{a(u^2+9u-3a)^3}{(u^3-2au^2-9au-2a^2-27a)^2}.
\]

\section{Field intersection problem: the case of $D_5$}\label{seQuin}
Let $\sigma:=(12345)$, $\rho:=(1243)$, $\tau:=\rho^2$, $\omega:=(12)\in S_5$ act on 
$k(x_1,\ldots,x_5)$ by $\pi(x_i)=x_{\pi(i)}, (\pi \in S_5)$. 
For simplicity, we write 
\[
C_5=\langle\sigma\rangle,\quad D_5=\langle\sigma,\tau\rangle,\quad 
F_{20}=\langle\sigma,\rho\rangle,\quad S_5=\langle\sigma,\omega\rangle.
\]
We take the cross-ratios 
\begin{align*}
x:=\xi(x_1,\ldots,x_5)=\frac{x_1-x_4}{x_1-x_3}\bigg{/}\frac{x_2-x_4}{x_2-x_3},\quad 
y:=\eta(x_1,\ldots,x_5)=\frac{x_2-x_5}{x_2-x_4}\bigg{/}\frac{x_3-x_5}{x_3-x_4}.
\end{align*}
Then $S_5$ acts faithfully on $k(x,y)$ as 
\begin{align*}
\sigma\ &:\ x\longmapsto y,\ y\longmapsto -(y-1)/x,&
\rho\ &:\ x\longmapsto x/(x-1),\ y\longmapsto (y-1)/(x+y-1),\\
\tau\ &:\ x\longmapsto x,\ y\longmapsto -(x-1)/y,& 
\omega\ &:\ x\longmapsto 1/x,\ y\longmapsto (x+y-1)/x.
\end{align*}
We take a $k$-generic polynomial $f_{s,t}^{D_5}(X)\in k(s,t)[X]$ 
as the formal $D_5$-relative $\langle\tau\rangle$-invariant resolvent polynomial by $x$: 
\begin{align*}
f_{s,t}^{D_5}(X)&:=\mathcal{RP}_{x,D_5}(X)
=\bigl(X-x\bigr)\bigl(X-y\bigr)\bigl(X-\frac{1-y}{x}\bigr)
\bigl(X-\frac{x+y-1}{xy}\bigr)\bigl(X-\frac{1-x}{y}\bigr)\\
&=X^5+(t-3)X^4+(s-t+3)X^3+(t^2-t-2s-1)X^2+sX+t.
\end{align*}
Note that $k(s,t)=k(x,y)^{D_5}$. 
We take two fields $k(\bx)=k(x,y)$ and $k(\bx')=k(x',y')$ where 
$x':=\xi(y_1,\ldots,y_5)$, $y':=\eta(y_1,\ldots,y_5)$ and the interchanging involution
\[
\iota \ :\ k(\bx,\bx')\longrightarrow k(\bx,\bx'),\quad
(x,y,x',y')\longmapsto (x',y',x,y)
\]
which is induced by the $\iota_{\bx,\by}$ of (\ref{defiota}).
We put $(s',t',d'):=\iota(s,t,d)$, $\bs':=(s',t')$ and $(\sigma',\tau',\rho')$ $:=$ 
$(\iota^{-1}\sigma\iota,\iota^{-1}\tau\iota$, $\iota^{-1}\rho\iota)$ $\in$ 
$\mathrm{Aut}_k(k(x',y'))$ and write 
\[
{D_5}'=\langle\sigma',\tau'\rangle,\quad 
{F_{20}}'=\langle\sigma',\rho'\rangle\quad \mathrm{and}\quad 
{D_5}''=\langle\sigma\sigma',\tau\tau'\rangle.
\]
We now take the $D_5\times D_5'$-primitive $D_5''$-invariant 
\begin{align*}
P:=xx'+yy'+\frac{(y-1)(y'-1)}{xx'}+\frac{(x+y-1)(x'+y'-1)}{xx'yy'}+\frac{(x-1)(x'-1)}{yy'}
\end{align*}
and the formal $D_5\times D_5'$-relative $D_5''$-invariant resolvent polynomial by $P$: 
\begin{align}
F_{\bs,\bs'}^1(X)&:=\mathcal{RP}_{P,D_5\times D_5'}(X)
=\prod_{\opi\in (D_5\times D_5')/D_5''}(X-\pi(P))\label{exform}\\
&\ =\Bigl(X^5-(t-3)({t'}-3)X^4+c_3X^3+\frac{c_2}{2}X^2+\frac{c_1}{2}X+\frac{c_0}{2}\Bigr)^2
\nonumber\\
&\qquad -\frac{d^2d'^2}{4}\Bigl(X^2+(t+{t'}-1)X+(s-t+{s'}-{t'}+t{t'}+2)\Bigr)^2;\nonumber
\end{align}
here $c_3,c_2,c_1,c_0\in k(s,t,{s'},{t'})$ are given by 
{\small 
\begin{align*}
c_3&=\bigl[2s-21t+3t^2-2t{s'}+t^2{\!}{s'}-t^2{t'}\bigr]+31-3s{s'}+5t{t'},\\
c_2&=\bigl[-20s+112t+8st-32t^2+2t^3+5t{s'}-13st{s'}-12t^2{\!}{s'}+4t^3{\!}{s'}-15st{t'}\\
&+14t^2{t'}+2t^3{t'}+8t^2{\!}{s'}{\!}{t'}-2t^3{t'}^2\bigr]-102+27s{s'}-119t{t'}
-st{s'}{t'}+6t^2{t'}^2,\\
c_1&=\bigl[32s+2s^2-128t-26st+60t^2+4st^2-8t^3-6s^2{s'}-7t{s'}+38st{s'}+9t^2{\!}{s'}
-5st^2{\!}{s'}\\
&-12t^3{\!}{s'}+2t^4{\!}{s'}-20t{s'}^2-8st{s'}^2+6t^2{\!}{s'}^2+2t^3{\!}{s'}^2+2st{t'}-77t^2{t'}
+3st^2{t'}+8t^3{t'}-29t^2{\!}{s'}{t'}\\
&+st^2{\!}{s'}{t'}+18t^3{\!}{s'}{t'}-2st^2{t'}^2+10t^3{t'}^2\bigr]+80-37s{s'}+145t{t'}
-45st{s'}{t'}+24t^2{t'}^2-8t^3{t'}^3,\\
c_0&=\bigl[-16s-2s^2+56t+24st+2s^2t-38t^2-8st^2+8t^3+5s^2{s'}-2t{s'}-38st{s'}-7s^2t{s'}\\
&+5t^2{\!}{s'}+13st^2{\!}{s'}+8t^3{\!}{s'}+2st^3{\!}{s'}-4t^4{\!}{s'}-21t{s'}^2-11st{s'}^2
-2t^2{\!}{s'}^2+2st^2{\!}{s'}^2+4t^3{\!}{s'}^2\\
&-104st{t'}-33s^2t{t'}+105t^2{t'}+35st^2{t'}+4t^3{t'}+16st^3{t'}-6t^4{t'}-2t^5{t'}
-s^2t{s'}{t'}+36t^2{\!}{s'}{t'}\\
&-14st^2{\!}{s'}{t'}-6t^3{\!}{s'}{t'}+6t^4{\!}{s'}{t'}+8t^2{\!}{s'}^2{t'}-37st^2{t'}^2
+22t^3{t'}^2-2st^3{t'}^2+8t^4{t'}^2+8t^3{\!}{s'}{t'}^2\\
&-2t^4{t'}^3\bigr]-24+14s{s'}-8s^2{s'}^2-224t{t'}+st{s'}{t'}
-101t^2{t'}^2-st^2{s'}{t'}^2-8t^3{t'}^3
\end{align*}
}
\noindent 
where $\bigl[a\bigr]:=a+\iota(a)$ for $a\in k(s,t,s',t')$, and 
$d^2\in k(s,t)$ is given by the formula 
\[
d^2=s^2-4s^3+4t-14st-30s^2t-91t^2-34st^2+s^2t^2+40t^3+24st^3+4t^4-4t^5. 
\]

We also take $\rho(P)\in k(s,t)$ which is conjugate of $P$ under the action of 
$F_{20}\times F_{20}'$ but not so under the action of $D_5\times D_5'$. 
Put 
\begin{align*}
F_{\bs,\bs'}^2(X)&:=\mathcal{RP}_{\rho(P),D_5\times D_5'}(X)
=F_{\rho(s),\rho(t),s',t'}^1(X)=\rho\bigl(F_{\bs,\bs'}^1(X)\bigr). 
\end{align*}

For $\ba=(a_1,a_2)\in M^2$, we denote $L_\ba:=\mathrm{Spl}_M f_\ba^{D_5}(X)$, 
$G_\ba:=\mathrm{Gal}(f_\ba^{D_5})$ and $G_{\ba,\ba'}:=\mathrm{Gal}(f_\ba^{D_5}f_{\ba'}^{D_5})$. 
We now state the result of the dihedral quintic case. 
\begin{theorem}[\cite{HM-3}]\label{thD5}
For $\ba=(a_1,a_2)$, $\ba'=(a_1',a_2')\in M^2$, we assume $\#G_\ba\geq \#G_{\ba'}>1$. 
An answer to the field intersection problem for $f_{s,t}^{D_5}(X)$ 
is given by the decomposition types {\rm DT}$(F_{\ba,\ba'}^1)$ and {\rm DT}$(F_{\ba,\ba'}^2)$ 
over $M$ as Table $2$ shows. 
\end{theorem}
\begin{center}
{\small
{\rm Table} $2$\vspace*{2mm}\\
}
{\footnotesize 
\begin{tabular}{|c|c|c|l|l|l|}\hline
$G_\ba$& $G_{\ba'}$ & $G_{\ba,{\ba'}}$ & & ${\rm DT}(F_{\ba,\ba'}^1)$ 
& ${\rm DT}(F_{\ba,\ba'}^2)$ \\ \hline 
& & $D_5\times D_5$ & $L_\ba\cap L_{\ba'}=M$ & $10$ & $10$ \\ \cline{3-6} 
& \raisebox{-1.6ex}[0cm][0cm]{$D_5$}  & $(C_5\times C_5)\rtimes C_2$ & 
$[L_\ba\cap L_{\ba'} : M]=2$ & $5,5$ & $5,5$\\ \cline{3-6}
& & \raisebox{-1.6ex}[0cm][0cm]{$D_5$} & \raisebox{-1.6ex}[0cm][0cm]{$L_\ba=L_{\ba'}$} 
& $5,2,2,1$ & $5,5$ \\ \cline{5-6}
$D_5$ & & & & $5,5$ & $5,2,2,1$ \\ \cline{2-6}
& $C_5$ & $D_5\times C_5$ & $L_\ba\cap L_{\ba'}=M$ & $10$ & $10$ \\ \cline{2-6}
& \raisebox{-1.6ex}[0cm][0cm]{$C_2$} & $D_{10}$ 
& $L_\ba\not\supset L_{\ba'}$ & $10$ & $10$ \\ \cline{3-6} 
& & $D_5$ & $L_\ba\supset L_{\ba'}$ & $5,5$ & $5,5$ \\ \cline{1-6}
& & $C_5\times C_5$ & $L_\ba\neq L_{\ba'}$ & $5,5$ & $5,5$\\ \cline{3-6}
\raisebox{-1.6ex}[0cm][0cm]{$C_5$} & $C_5$  & \raisebox{-1.6ex}[0cm][0cm]{$C_5$} 
& \raisebox{-1.6ex}[0cm][0cm]{$L_\ba=L_{\ba'}$} 
& $5,1,1,1,1,1$ & $5,5$ \\ \cline{5-6}
& & & & $5,5$ & $5,1,1,1,1,1$ \\ \cline{2-6}
& $C_2$ & $C_{10}$ & $L_\ba\cap L_{\ba'}=M$ & $10$ & $10$ \\ \hline
\raisebox{-1.6ex}[0cm][0cm]{$C_2$} & \raisebox{-1.6ex}[0cm][0cm]{$C_2$} & $C_2\times C_2$ 
& $L_\ba\neq L_{\ba'}$ & $4,4,2$ & $4,4,2$\\ \cline{3-6}
& & $C_2$ & $L_\ba=L_{\ba'}$ & $2,2,2,2,1,1$ & $2,2,2,2,1,1$\\ \hline
\end{tabular}
}\vspace*{0mm}
\end{center}
We checked the decomposition types on Table $2$ using the computer algebra system 
GAP \cite{GAP}.
\begin{remark}
In the cases of $G_\ba=D_5$, $C_2$, 
the quadratic subextension of $L_{s_1,t_1}$ over $M$ for $(s_1,t_1)\in M^2$ is given by 
{\small 
\begin{align*}
M\Bigl(\sqrt{s_1^2-4s_1^3+4t_1-14s_1t_1-30s_1^2t_1-91t_1^2-34s_1t_1^2+s_1^2t_1^2
+40t_1^3+24s_1t_1^3+4t_1^4-4t_1^5}\Bigr). 
\end{align*}
}
\end{remark}
\begin{remark}
As is well known, for a suitable integer $l$, we may have an answer to the field intersection 
problem via the formal $D_5\times D_5'$-relative 
$\langle\tau\rangle\times\langle\tau'\rangle$-invariant resolvent polynomial 
$G_{\bs,\bs'}(X):=\mathcal{RP}_{x+lx'\!,D_5\times D_5'}(X)$ by $x+lx'$ 
(cf. \cite{Coh93}, \cite{Coh00}). 
Although the degree of the multi-resolvent polynomials $F_{\bs,\bs'}^1(X)$ and 
$F_{\bs,\bs'}^1(X)$ with respect to $X$ is $10$, the degree of $G_{\bs,\bs'}(X)$ is $25$ and 
an explicit formula of $G_{\bs,\bs'}(X)$ in terms of $\bs$ and $\bs'$ is very complicated. 
We remark that to construct a suitable explicit formula is significant to investigate the 
structure of all $G$-Galois extensions over $M$ (cf. Section \ref{seIntro} and Corollaries 
\ref{cor0} and \ref{cor01}). 
\end{remark}
\begin{example}
Take $M=\mathbb{Q}$, $\ba=(0,3)$ and $\ba'=(10,3)$. 
Then we have 
{\small 
\begin{align*}
F_{\ba,\ba'}^1(X)&=(X+5)^3(X^2-15X+150)(X^5-625X^2-9375),\\
F_{\ba,\ba'}^2(X)&=\Bigl(X^5-\frac{125}{9}X^3+\frac{6250}{81}X-\frac{3125}{27}\Bigr)
\Bigl(X^5-\frac{125}{9}X^3+\frac{625}{9}X^2+\frac{15625}{81}X+\frac{15625}{27}\Bigr).
\end{align*}
}
\hspace*{-3.2mm} 
The decomposition type of $F_{\ba,\ba'}^1(X)$ over $\mathbb{Q}$ should be $5,2,2,1$ 
(cf. also Theorem \ref{thfun}); 
and hence we conclude that $\mathrm{Spl}_\mathbb{Q} f_\ba^{D_5}(X)=\mathrm{Spl}_\mathbb{Q} 
f_{\ba'}^{D_5}(X)$ and $G_\ba=D_5$. 
\end{example}
From the viewpoint of Diophantine geometry, we give some numerical examples 
of the field isomorphism problem of $f_{s,t}^{D_5}(X)$ over $M=\mathbb{Q}$ 
and for integral points $\ba,\ba'\in\mathbb{Z}^2$ using Theorem \ref{thD5} and 
the explicit formula (\ref{exform}). 
We do not know, however, for a given $\ba\in\mathbb{Z}^2$ 
whether there exist only finitely many $\ba'\in\mathbb{Z}^2$ such that 
$\mathrm{Spl}_\mathbb{Q} f_\ba^{D_5}(X)=\mathrm{Spl}_\mathbb{Q} f_{\ba'}^{D_5}(X)$ or not 
(cf. Corollary \ref{cor01}). 

\begin{example}
Take $M=\mathbb{Q}$ and $t:=1$. 
Then we have $f_{s,1}^{D_5}(X)=X^5-2X^4+(s+2)X^3-(2s+1)X^2+sX+1$. 
For $s_1,s_1'\in \mathbb{Z}$ in the range $-10000\leq s_1< s_1'\leq 10000$, 
$\mathrm{Spl}_\mathbb{Q} f_{s_1,1}^{D_5}(X)=\mathrm{Spl}_\mathbb{Q} f_{s_1',1}^{D_5}(X)$ if and 
only if $(s_1,s_1')\in X_1\cup X_2$ where
\begin{align*}
\vspace*{-1mm}
X_1&=\{(-6,0),(-1,41),(-94,-10)\},\\
X_2&=\{(-1,0),(-6,-1),(-18,-7),(1,34),(0,41),(-6,41),(-167,-8)\}.
\end{align*}
It was directly checked by Theorem \ref{thD5} that, in the range $-10000\leq s_1<s_1'\leq 10000$, 
$(s_1,s_1')\in X_i$ if and only if $\mathrm{DT}(F_{s_1,1,s_1',1}^i/\mathbb{Q})$ 
includes $1$, for each of $i=1,2$. 
\end{example}
\begin{example}
Kida-Renault-Yokoyama \cite{KRY} showed that there exist infinitely many $b\in \mathbb{Q}$ 
such that $\mathrm{Spl}_\mathbb{Q} f_{0,1}^{D_5}(X)=\mathrm{Spl}_\mathbb{Q} f_{b,1}^{D_5}(X)$. 
Their method enables us to construct such $b$'s explicitly via rational points of an 
associated elliptic curve (cf. \cite{KRY}). 
They also pointed out that in the range $-400\leq s_1, t_1\leq 400$ there are $25$ pairs 
$(s_1,t_1)\in \mathbb{Z}^2$ such that $\mathrm{Spl}_\mathbb{Q} f_{0,1}^{D_5}(X)
=\mathrm{Spl}_\mathbb{Q} f_{s_1,t_1}^{D_5}(X)$. 
We may classify the $25$ pairs by the polynomials $F_{0,1,s_1,t_1}^1(X)$ and 
$F_{0,1,s_1,t_1}^2(X)$. 
In the range above, for $i=1,2$, 
$\mathrm{DT}(F_{0,1,s_1,t_1}^i/\mathbb{Q})$ includes $1$ if and only if 
$(s_1,t_1)\in X_i$ where 
\vspace*{-1mm}
\begin{align*}
X_1=\{&(0,1),(4,-1),(4,5),(-6,1),(-24,19),(34,11),(36,-5),\\
&(46,-1),(-188,23),(264,31),(372,-5),(378,43)\},\\
X_2=\{&(-1,-1),(-1,1),(5,-1),(41,1),(-43,5),(47,13),(59,-5),\\
&(59,19),(101,19),(125,-23),(149,11),(155,25),(-169,55)\}.
\end{align*}
\vspace*{-6mm}\\
By Theorem \ref{thD5}, we checked such pairs in the range $-20000\leq s_1, t_1\leq 20000$; 
and added six pairs $(526,41)$, $(952,113)$, $(2302,95)$, $(6466,311)$, $(7180,143)$ and 
$(7480,-169)$ to $X_1$ and just four pairs $(785,-25)$, $(3881,29)$, $(-11215,299)$ and 
$(19739,-281)$ to $X_2$. 
\end{example}
\begin{example}
Take the $k$-generic polynomial $g_{A,B}^{C_5}(X)$ $\in$ $k(A,B)[X]$ for $C_5$ 
which is constructed by Hashimoto-Tsunogai \cite{HT03}: 
\vspace*{-1mm}
\begin{align*}
g_{A,B}^{C_5}(X)=X^5-\frac{P}{Q^2}(A^2-2A+15B^2+2)X^3+\frac{P^2}{Q^3}(2BX^2-(A-1)X-2B) 
\end{align*}
\vspace*{-5mm}\\
where $P=(A^2-A-1)^2+25(A^2+1)B^2+125B^4$, $Q=(A+7)B^2-A+1$. 
We may apply Theorem \ref{thD5} to $g_{A,B}^{C_5}(X)$ since 
there exist $s_1,t_1\in k(A,B)$ such that 
$\mathrm{Spl}_{k(A,B)} f_{s_1,t_1}^{D_5}(X)=\mathrm{Spl}_{k(A,B)} g_{A,B}^{C_5}(X)$ 
(cf. \cite{HT03}, \cite{HM-3}). 
For $\ba=(a,b)$, $\ba'=(a',b')\in \mathbb{Z}^2$, 
if $\ba'=(a,\pm b)$ or $\{\ba,\ba'\}=\{(-1,\pm b),(1,\pm b)\}$ then 
$\mathrm{Spl}_\mathbb{Q} g_{a,b}^{C_5}(X)=\mathrm{Spl}_\mathbb{Q} g_{a',b'}^{C_5}(X)$ 
(cf. \cite{HM-3}). 
For $\ba, \ba'\in \mathbb{Z}^2$ in the range $-50\leq a, a'\leq 50$, $0\leq b\leq b'\leq 50$ 
with $\ba\neq \ba'$, $\{\ba,\ba'\}\neq\{(-1,b),(1,b)\}$, we see 
that $\mathrm{Spl}_\mathbb{Q} g_{a,b}^{C_5}(X)=\mathrm{Spl}_\mathbb{Q} g_{a',b'}^{C_5}(X)$ 
if and only if $(a,b,a',b')\in X_1\cup X_2$ where
\vspace*{-1mm}
\begin{align*}
X_1&=\{(-3,1,-3,11),(3,3,23,3),(23,3,3,3),(7,3,27,9),\\
&\qquad (2,2,-28,14),(8,11,33,14),(23,5,35,7),(41,11,-15,17)\},\\
X_2&=\{(-2,1,3,2),(4,1,-6,2),(3,1,13,7),(16,2,-12,5),(-2,2,18,4),(31,1,-19,7),\\
&\qquad (-3,3,-33,3),(-33,3,-3,3),(-16,13,34,19),(-2,3,43,6),(12,4,46,10)\}.
\end{align*}
By Theorem \ref{thD5}, it can be checked, in the range above and for each of $i=1,2$, that 
$(a,b,a',b')\in X_i$ if and only if the decomposition type of $F_{\ba,\ba'}^i(X)$ over 
$\mathbb{Q}$ includes $1$. 
\end{example}
\begin{example}
Let $h_n(X)$ be Lehmer's simplest quintic polynomial 
\vspace*{-1mm}
\begin{align*}
h_n(X)=X^5&+n^2X^4-(2n^3+6n^2+10n+10)X^3\\
&+(n^4+5n^3+11n^2+15n+5)X^2+(n^3+4n^2+10n+10)X+1
\end{align*}
\vspace*{-6mm}\\
(cf. \cite{Leh88}), and take $M=\mathbb{Q}$. We regard $n$ as an 
independent parameter over $\mathbb{Q}$. 
By the result in \cite{HR}, for Brumer's quintic $f_{s,t}^{D_5}(X)$, we see that 
$\mathrm{Spl}_{\mathbb{Q}(n)} h_n(X)=\mathrm{Spl}_{\mathbb{Q}(n)}f_{s,t}^{D_5}(X)$ 
where $s=-20-5n+10n^2+12n^3+5n^4+n^5$, $t=-7-10n-5n^2-n^3$. 
By Theorem \ref{thD5}, we checked pairs $(n, n')\in\mathbb{Z}^2$ in the range 
$-10000\leq n<n'\leq 10000$ to confirm that $\mathrm{Spl}_\mathbb{Q} h_n(X)
=\mathrm{Spl}_\mathbb{Q} h_{n'}(X)$ if and only if $(n,n')=(-2,-1)$. 
\end{example}
\begin{acknowledgments}
We thank Professor Kazuhiro Yokoyama for drawing our attention to multi-resolvent 
polynomials and also for his useful suggestions. 
We also thank the referee for giving fruitful comments and for careful reading of the 
manuscript. 
\end{acknowledgments}


\end{document}